\documentclass[11pt,leqno]{article}
\usepackage{amssymb,amsfonts}
\usepackage{amsmath,amsthm,amsxtra}
\usepackage{epsfig}
\usepackage{slashed}
\usepackage{color}
\usepackage{verbatim}

\newcommand{\ad}{\mathop{\rm ad}\,}

\newcommand{\ch}{{\rm ch}}

\newcommand{\End}{\mathop{\rm End }}

\renewcommand{\Im}{\mathop{\rm Im  \, }}

\newcommand{\re}{\mathop{\rm re  \, }}

\newcommand{\sdim}{\mathop{\rm sdim \, }}

\newcommand{\tr}{\mathrm{tr} \, }

\newcommand{\CC}{\mathbb{C}}

\newcommand{\QQ}{\mathbb{Q}}

\newcommand{\ZZ}{\mathbb{Z}}

\newcommand{\fg}{\mathfrak{g}}
\newcommand{\fh}{\mathfrak{h}}

\newcommand\bl{(\, . \, | \, . \, )}

\makeatletter
\renewcommand\section{\@startsection {section}{1}{\z@}%
                                 {-3.5ex \@plus -1ex \@minus -.2ex}%
                                   {2.3ex \@plus.2ex}%
                                   {\normalfont\large\bfseries}}
\renewcommand\subsection{\@startsection{subsection}{2}{\z@}%
                                     {-3.25ex\@plus -1ex \@minus -.2ex}%
                                     {0ex \@plus .0ex}%
                                     {\normalfont\normalsize\bfseries}}

\@addtoreset{equation}{section}
\makeatother

\newtheorem{theorem}{Theorem}

\newtheorem{proposition}[theorem]{Proposition}

\newtheorem*{lemma*}{Lemma}

\theoremstyle{remark}
\newtheorem{remark}{Remark}
\newtheorem{example}{Example}

\makeatletter
\def\@maketitle{\newpage
 \null
 \vskip 2em
 \begin{center}%
 \vskip 3em
  {\Large\bf \@title \par}%
  \vskip 1.5em
  {\normalsize
   \lineskip .5em
   \begin{tabular}[t]{c}\@author
   \end{tabular}\par}%
  \vskip 2em

 \end{center}%
 \par
 \vskip 2.5em}
\makeatother

\renewcommand{\epsilon}{\varepsilon}

\definecolor{light}{gray}{.9}

\newcommand{\half}{\frac{1}{2}}

\newcommand{\la}{\lambda}
\newcommand{\La}{\Lambda}
\newcommand{\al}{\alpha}

\newcommand{\wg}{\widehat{\fg}}
\newcommand{\wh}{\widehat{\fh}}
\newcommand{\wrh}{\widehat{\rho}}

\newcommand{\tz}{(\tau, z)}
\newcommand{\tzt}{(\tau, z, t)}

\newcommand{\LLa}{L(\Lambda)}
\newcommand{\tot}{\frac{\tau}{2}}

\begin{document}

\title{A remark on boundary level admissible representations}


\author{Victor G.~Kac\thanks{Department of Mathematics, M.I.T, 
Cambridge, MA 02139, USA. Email:  kac@math.mit.edu ~~~} \ 
and Minoru Wakimoto\thanks{Email: ~~wakimoto@r6.dion.ne.jp~~~}}

\maketitle



Recently a remarkable map between 4-dimensional superconformal field theories and vertex algebras has been constructed \cite{BLLPRV15}. This has lead to new insights in the theory of characters of vertex algebras. In particular it was observed that in some cases these characters decompose in nice products
\cite{XYY16}, \cite{Y16}. 

The purpose of this note is to explain the latter phenomena. Namely, we point out that it is immediate by our character formula \cite{KW88}, \cite{KW89} that in the case of a \textit{boundary level} the characters of admissible representations of affine Kac-Moody algebras and the corresponding $W$-algebras
decompose in products in terms of the Jacobi form $ \vartheta_{11} \tz. $

We would like to thank Wenbin Yan for drawing our attention to this question. 

Let $ \fg $ be a simple finite-dimensional Lie algebra over $ \CC, $ let $ \fh $ be a Cartan subalgebra of $ \fg,  $ and let $ \Delta \subset \fh^* $ be the set of roots. Let $ Q = \ZZ \Delta $ be the root lattice and let $ Q^* = \{ h \in \fh \ | \ \al (h) \in \ZZ \mbox{ for all } \al \in \Delta \}  $ be the dual lattice. 
Let $ \Delta_+\subset \Delta $ be a subset of positive roots, let $ \{ \al_1, \ldots, \al_\ell \} $ be the set of simple roots and let $  \rho  $  be half of the sum of positive roots. Let $ W $ be the Weyl group. Let $ \bl $ be the invariant symmetric bilinear form on $ \fg,  $ normalized by the condition $ (\al|\al) = 2 $ for a long root $ \al, $ and let $ h^\vee $ be the dual Coxeter number ($ = \half $ eigenvalue of the Casimir operator on $ \fg $). We shall identify $ \fh  $ with $ \fh^* $ using the form $ \bl. $

Let $ \wg = \fg [t, t^{-1}] + \CC K + \CC d $ be the associated to $ \fg $ affine Kac-Moody algebra (see \cite{K90} for details), let $ \wh = \fh + \CC K + \CC d $ be its Cartan subalgebra. We extend the symmetric bilinear form $ \bl $ from $ \fh  $ to $ \wh $ by letting $ (\fh | \CC K + \CC d) =0,  (K|K) = 0, (d|d) = 0, (d|K)= 1, $ and we identify $ \wh^*  $ with $ \wh $ using this form. Then $ d $ is identified with the $ 0^{th} $ fundamental weight $ \La_0 \in \wh^*, $ such that $ \La_0 |_{\fg [t, t^{-1}]+ \CC d} =0, \La_0 (K) = 1,  $ and $ K $ is identified with the imaginary root $ \delta \in \wh^*.  $
Then the set of real roots of $ \wg $ is $ \hat{\Delta}^{\re} = \{ \al + n \delta |\, \al \in \Delta, n \in \ZZ \}$ and the subset of positive real roots is $ \hat{\Delta}^{\re}_+ = \Delta_+ \cup \{\al + n \delta |\, \al \in \Delta, n \in \ZZ_{\geq 1} \} $. Let $\hat{\rho}=h^\vee \Lambda_0 +\rho$. Let
\[ \hat{\Pi}_u = \{ u \delta - \theta, \al_1, \ldots, \al_\ell \}, \]
where $ \theta \in \Delta_+ $ is the highest root, so that $ \hat{\Pi}_1 $ is the set of simple roots of $ \wg. $ For $ \al \in \hat{\Delta}^{\re} $ one lets $ \al^\vee = 2 \al / (\al|\al).  $ Finally, for $ \beta \in Q^* $ define the translation $ t_\beta \in \End \wh^* $ by 
\[ t_\beta (\la) = \la + \la (K) \beta -
((\la |\beta) + \half \la (K) |\beta|^2) \delta . 
 \]

 Given $ \La \in \wh^* $ let $ \hat{\Delta}^\La = \{ \al \in \hat{\Delta}^{\re}
 |\, (\La|\al^\vee) \in \ZZ  \}  $. Then $ \La $ is called an \textit{admissible} weight if the following two properties hold 
\begin{enumerate}
\item[(i)]  $ (\La + \wrh | \al^\vee) \notin \ZZ_{\leq 0}$ for all $ \al \in \hat{\Delta}_+, $
\item[(ii)] $ \QQ \hat{\Delta}^\La = \QQ \hat{\Delta}. $
\end{enumerate}
If instead of (ii) a stronger condition holds: 
\begin{enumerate}
\item[(ii)$ ' $] $ \varphi (\hat{\Delta}^\La) = \hat{\Delta} $  for a linear
  isomorphism $ \varphi : \wh^* \rightarrow \wh^*, $ 
\end{enumerate}
then $ \La  $ is called a \textit{principal} admissible weight. In \cite{KW89} the classification and character formulas for admissible weights is reduced to that for principal admissible weights. The latter are described by the following proposition. 

\begin{proposition}
\label{prop1} 
\cite{KW89}  Let $ \La $ be a principal admissible weight and let $ k = \La (K) $ be its level. Then
\begin{enumerate}
\item[(a)] $ k$ is a rational number with denominator $ u \in \ZZ_{\geq 1}, $ such that
\begin{equation}
\label{1} 
k + h^\vee \geq \frac{h^\vee}{u} \mbox{ and } \gcd (u, h^\vee) = \gcd (u, r^\vee) = 1, 
\end{equation}
where $ r^\vee = 1 $ for $ \fg  $ of type A-D-E, = 2 for $ \fg  $ of type B, C, F, and = 3 for $ \fg = G_2. $
\item[(b)] All principal admissible weights are of the form
\begin{equation} 
\label{2}
\La = (t_\beta y). (\La^0 - (u-1) (k+h^\vee)\La_0), 
\end{equation}
where $ \beta \in Q^*, y \in W $ are such that $ (t_\beta y) \hat{\Pi}_u  \subset \hat{\Delta}_+,  \La^0$ is an integrable weight of level $ u(k+h^\vee)-h^\vee, $ and dot denotes the shifted action: $ w.\La = w(\La + \wrh) - \wrh. $
\item[(c)] For $ \fg = s\ell_N $ all admissible weights are principal admissible. 
\end{enumerate}
\end{proposition}

Recall that the normalized character of an irreducible highest weight $ \wg $-module $ \LLa $ of level $ k \neq -h^\vee $ is defined by
\[ \ch_\La (\tau, z, t) = q^{m_\La} \tr_{\LLa} e^{2 \pi i h} \]
where 
\begin{equation} 
\label{3}
 h = -\tau d + z + tK, \  z \in \fh, \ \tau, t \in \CC, \ \Im \tau > 0, \ q = e^{2 \pi i \tau},  
\end{equation} 
and $ m_\La = \frac{|\La + \wrh |^2}{2 (k+h^\vee)} -\frac{\dim \fg}{24} $ (the normalization factor $ q^{m_\La} $ ``improves'' the modular invariance of the character).

In \cite{KW89} the characters of the $ \wg $-modules $ \LLa $ for arbitrary admissible $ \La  $ were computed, see Theorem 3.1, or formula (3.3) there for another version in case of a principal admissible $ \La.  $
In order to write down the latter formula, recall the normalized affine denominator for $ \wg: $
\[ \hat{R} (h) = q^{\frac{\dim \fg}{24}} e^{\wrh(h)} \prod_{n=1}^{\infty} (1-q^n)^\ell \prod_{\al \in \Delta_+} (1-e^{\al(z)}q^n) (1-e^{-\al(z)} q^{n-1}). \]
In coordinates \eqref{3} this becomes:
\begin{equation}
\label{4}
\hat{R} \tzt = (-i)^{|\Delta_+|} e^{2 \pi i h^\vee t} \eta (\tau)
^{\half(3\ell - \dim \fg)} 
\prod_{\al \in \Delta_+} \vartheta_{11} (\tau, \al(z)),
\end{equation}
where 
\[ \vartheta_{11} \tz =- i q^{\frac{1}{12}} e^{-\pi i z} \eta (\tau) \prod_{n=1}^{\infty} (1-e^{-2 \pi i z} q^n)(1-e^{2 \pi i z} q^{n-1}) \]
is one of the standard Jacobi forms  $\vartheta_{ab}$, $a,b=0$ or 1 
(see e.g., Appendix to \cite{KW14}), and $\eta(\tau)$ is the Dedekind eta function.

For a principal admissible $ \La,  $ given by \eqref{2}, formula (3.3) from \cite{KW89} becomes in coordinates \eqref{3}:
\begin{equation}
\label{5}
(\hat{R} \ch_\La) \tzt = (\hat{R} \ch_{\La^0}) \left( u \tau, y^{-1} (z + \tau \beta), \frac{1}{u} (t + (z|\beta) + \frac{\tau |\beta|^2}{2}) \right).
\end{equation}
It follows from \eqref{5} that if $ \La^0 = 0 $ in \eqref{2} (so that $ \ch_{\La^0} =1 $), which is equivalent to
\begin{equation}
\label{6}
k + h^\vee = \frac{h^\vee}{u} \mbox{ and } \gcd (u, h^\vee) = \gcd (u, r^\vee) = 1,
\end{equation}
the (normalized) character $ \ch_\La $ turns into a product. The level $ k, $ defined by \eqref{6}, is naturally called the \textit{boundary principal admissible} level in \cite{KRW03}, see formula (3.5) there. We obtain from Proposition \ref{prop1}, \eqref{4} and \eqref{5}
\begin{proposition}
\label{prop2}
\begin{enumerate}
\item[(a)]  All boundary principal admissible weights are of level $ k,  $ given by \eqref{6}, and are of the form
\begin{equation}
\label{7}
\La = (t_\beta y). (k \La_0),
\end{equation}
where $ \beta \in Q^*, y \in W $ are such that $ (t_\beta y) \hat{\Pi}_u \subset \hat{\Delta}_+. $ In particular, $ k \La_0 $ is a principal admissible weight of level \eqref{6}.
\item[(b)] If $ \La $  is of the form \eqref{7}, then
\[ \ch_\La \tzt = e^{2 \pi i (kt + \frac{h^\vee}{u}(z | \beta))} q^{\frac{h^\vee}{2u}|\beta|^2} \left( \frac{\eta (u \tau)}{\eta (\tau)}\right)
^{\half(3\ell - \dim \fg)}  \prod_{\al \in \Delta_+} \frac{\vartheta_{11} (u \tau, y(\al) (z + \tau \beta))}{\vartheta_{11} (\tau, \al (z))}.  \]
\end{enumerate}
\end{proposition}

\begin{remark}
\label{rem1}
For the vacuum module $ L(k\La_0) $ of the boundary principal admissible 
level $k$ the character formula from Proposition \ref{prop2}(b) becomes
\[ \ch_{k\La_0} \tzt = e^{2 \pi i kt} \left( \frac{\eta(u \tau)}{\eta (\tau)}\right)^{\half(3\ell - \dim \fg)}  \prod_{\al \in \Delta_+} \frac{\vartheta_{11} (u \tau, \al(z))}{\vartheta_{11} (\tau, \al(z))}. \]
\end{remark}

\begin{example}
\label{ex1}
Let $ \fg = s \ell_2, $ so that $ h^\vee = 2.  $ Then the boundary levels are $ k = \frac{2}{u}-2, $ where $ u $ is a positive odd integer, and all admissible weights are 
\[  \La_{k,j} : = t_{-\frac{j}{2} \al_1} . (k \La_0) = (k + \frac{2j}{u}) \La_0 - \frac{2j}{u} \La_1, \ j = 0,1, \ldots, u-1,  \]
and the character formula from Proposition \ref{prop2}(b) becomes:
\begin{equation}
\label{8}
\ch_{\La_{u,j}} = e^{2 \pi i (kt - \frac{j}{u}z)} q^{\frac{j^2}{2u}} \frac{\vartheta_{11} (u \tau, z-j\tau)}{\vartheta_{11} \tz}.
\end{equation}
For $ u = 3 $ and 5 some of these formulas were conjectured in \cite{Y16}.
\end{example}

\begin{example}
\label{ex2}
Let $ \fg = s\ell_N, $ so that $ h^\vee = N,  $ let $ N>1 $ be odd, and let $ u =2 $. Then the boundary admissible level is $ k = -\frac{N}{2}, $ and the boundary admissible weights of the form $ t_\beta . (k\La_0) $ are:
\[ \La_{N,p} = -\frac{N}{2} \La_p, \ p = 0,1, \ldots, ,N-1, \]
where $ \La_p $ are the fundamental weights of $ \wg. $
Letting $ z = \sum_{i =1}^{N-1} z_i \bar{\La}_i, $ where $ \bar{\La}_i $ are the fundamental weights of $ \fg, $ the character formula from Proposition \ref{prop2} (b) becomes:
\[ 
\begin{aligned}
&\ch_{\La_{N,p}} \tzt = i^{p(N-p)} e^{-\pi i Nt} \left( \frac{\eta (2 \tau)}{\eta (\tau)}\right)^{- \frac{(N-1)(N-2)}{2}} \\
&\times \frac{ \displaystyle \prod_{\substack{ 1 \leq i \leq j < p \\ \mbox{ or } p < i \leq j < N}} \vartheta_{11} (2 \tau, z_i + \ldots + z_j) \prod_{1 \leq i \leq p  \leq j < N} \vartheta_{01} (2 \tau, z_i + \ldots + z_j)}{\displaystyle \prod_{1 \leq i  \leq j < N} \vartheta_{11} ( \tau, z_i + \ldots + z_j)}, 
\end{aligned} \]
where 
\[ \vartheta_{01} \tz = \prod_{n =1 }^{\infty} (1-q^n) (1-e^{2 \pi i z} q^{n-\half}) (1-e^{-2 \pi i z} q^{n-\half}). \]
This follows from Proposition \ref{prop2}(b) by applying to $ \vartheta_{11} $ an elliptic transformation (see e.g. \cite{KW14}, Appendix).
In particular
\[ \ch_{-\frac{N}{2} \La_0} = e^{- \pi i Nt} \left(\frac{\eta (2 \tau)}{\eta (\tau)}\right)^{-\frac{(N-1)(N-2)}{2}} \prod_{1 \leq i \leq j < N} \frac{\vartheta_{11} (2 \tau, z_i + \ldots + z_j)}{\vartheta_{11} (\tau, z_i+ \ldots + z_j)}. \]
The latter formula was conjectured in \cite{XYY16}. 
\end{example}

\begin{remark}
\label{rem2}
For principal admissible weights $ \La = (t_\beta y). (k\La_0) $ and $ (t_{\beta'} y'). (k\La_0) $ of boundary level $ k = \frac{h^\vee}{u}-h^\vee  $ the $ S $-transformation matrix $ (a(\La, \La')), $ given by \cite{KW89}, Theorem 3.6, simplifies to  
\[ a(\La, \La') = | Q / uh^\vee Q^* |^{-\half} \epsilon(yy') \prod_{\al \in \Delta_+} 2 \sin \frac{\pi iu  (\rho | \al)}{h^\vee} e^{-2\pi i \left((\rho|\beta + \beta') + \frac{h^\vee (\beta| \beta')}{u}\right)} . \]
\end{remark}

\begin{remark}
\label{rem3}
If $ \fg = s \ell_2 $ and $ k $ is as in Example \ref{ex1}, then 
\[ a (\La_{k,j}, \La_{k,j'}) = (-1)^{j + j'} e^{- \frac{2 \pi i jj'}{u}} \frac{1}{\sqrt{u}} \sin \frac{u \pi }{2}. \]
One can compute fusion coefficients by Verlinde's formula:
\[ N_{\La_{k, j_1}, \La_{k, j_2}, \La_{k, j_3}} = (-1)^{j_1+j_2+j_3} \mbox{ if } j_1+j_2+j_3 \in u \ZZ, \mbox{ and } = 0 \mbox{ otherwise}. \]
\end{remark}

\begin{example} 
\label{ex3}
Let $\fg=sl_3$, so that $h^\vee=3$, and let $u$ be a positive integer, coprime to 3. Then all (principal) admissible weights have level $k=\frac{3}{u}-3$ 
and are of the form \eqref{7}, where
\[\beta=-(-1)^p
(k_1\bar{\Lambda}_1+k_2\bar{\Lambda}_2),\, y=r_{\theta}^p,\, p=0\, 
\mbox{or}\,1, \,k_i\in\ZZ, k_i\geq \delta_{p,1},\, k_1+k_2\leq u-\delta_{p,0}.\]
Denote this weight by $\Lambda^{(p)}_{u;k_1,k_2}=(t_\beta y).(k\Lambda_0)$. Using Remark \ref{rem2}, one computes the fusion coefficients by Verlinde's formula:
\[N_{\Lambda^{(p)}_{u;k_1,k_2}\Lambda^{(p')}_{u;k'_1,k'_2}\Lambda^{(p'')}_
{u;k''_1,k''_2}}=(-1)^{p+p'+p''} \,\mbox{if}\,\,
(-1)^{p}k_i+(-1)^{p'}k'_i+(-1)^{p''}k''_i\in u\ZZ \,\,\mbox{for}\,i=1,2,\] 
and $=0$ otherwise.
\end{example} 

\begin{remark}
\label{rem 4}
If $ \La $ is an arbitrary admissible weight, then $ \hat{\Delta}^\La $ decomposes in a disjoint union of several affine root systems. Then $ \La $ has \textit{boundary level} if restrictions of it to each of them has boundary level, and formula (3.4) from \cite{KW89} shows that $ \ch_\La $ decomposes in a product of the corresponding boundary level characters. Note also that all the above holds also for twisted affine Kac-Moody algebras \cite{KW89}. 
\end{remark} 

\begin{remark}
\label{rem 5}
The product character formula for boundary level affine Kac-Moody superalgebras holds as well, see \cite{GK15}, formula (2). 
\end{remark}

Recall that to any $ s\ell_2$-triple $ \{ f, x, e \}  $ in $ \fg, $ where $ [x,f] = -f, \ [x,e] = e, $ one associates a $ W $-algebra $ W^k (g,f) $, obtained from the vacuum $ \wg $-module of level $ k $ by quantum Hamiltonian reduction, so that any $ \wg $-module $ \LLa $ of level $ k $ produces either an irreducible $ W^k (g,f) $-module $ H(\La) $ or zero. The characters of $ \LLa $ and $ H(\La) $ are related by the following simple formula (\cite{KRW03} or \cite{KW14}):
\begin{equation}
\label{9}
\left( \overset{W}{R} \ch_{H(\La)} \right) \tz = \left(\hat{R} \ch_\La \right) (\tau, -\tau x + z, \tot (x|x)).
\end{equation}
Here $ z \in \fh^f $, the centralizer of $f$ in $\fh$, and 
\begin{equation}
\label{10}
\overset{W}{R} \tz = \eta (\tau)^{\frac{3}{2}l - \half \dim (\fg_0 + \fg_{1/2})} \prod_{\al \in \Delta^0_+} \vartheta_{11} (\tau, \al (z)) \left(\prod_{\al \in \Delta_{1/2}}  \vartheta_{01} (\tau, \al(z)) \right)^{1/2},
\end{equation}
where $ \fg = \oplus_j \fg_j $ is the eigenspace decomposition for $ \ad x, \ \Delta_j \subset \Delta $ is the set of roots of root spaces in $ \fg_j $ and $ \Delta^0_+ = \Delta_+ \cap \Delta_0 $ (we assume that $ \Delta_j \subset \Delta_+ $ for $ j > 0 $). If $ k $ is a boundary level \eqref{6}, we obtain from Proposition \ref{prop2}(b) and formulas \eqref{9}, \eqref{10} the following character formula for $ H(\La) $ if $ \La $ is a principal admissible weight \eqref{7}  ($z\in \fh^f$):
\begin{equation}
\label{11}
\begin{aligned}
 \ch_{H(\La)} \tz & = (-i)^{|\Delta_+|}q^{\frac{h^\vee}{2u} |\beta - x|^2} e^{\frac{2 \pi i h^\vee}{u} (\beta | z)} \\
& \times \frac{\eta (u \tau)^{\frac{3}{2}\ell - \half \dim \fg}}{\eta (\tau)^{\frac{3}{2}\ell - \half \dim (\fg_0 + \fg_{1/2})}} \  \frac{\displaystyle \prod_{\al \in \Delta_+} \vartheta_{11} (u \tau, y (\al) (z + \tau \beta -\tau x))}{ \displaystyle \prod_{\al \in \Delta^0_+} \vartheta_{11} (\tau, \al(z)) \left( \displaystyle \prod_{\al \in \Delta_{1/2}} \vartheta_{01} (\tau, \al(z))\right)^{1/2}}.
\end{aligned}
\end{equation}

\begin{remark}
\label{rem 6}
A formula, similar to Proposition \ref{prop2}(b) and to formula \eqref{11}, holds if $ \fg $ is a basic Lie superalgebra; one has to replace the character by the supercharacter, $ \dim  $ by $ \sdim, $ and the factor $ \vartheta_{ab},  $ corresponding to a root $ \al,  $ by its inverse if this root is odd.
Also, the character is obtained from the supercharacter by replacing
$\vartheta_{ab}$  by $\vartheta_{a,b+1\!\! \mod 2}$ if the root $\alpha$ is odd.  
\end{remark}

\begin{remark}
\label{rem 7}
An example of \eqref{11} is the minimal series representations of the Virasoro algebra with
central charge $c=1-\frac{3(u-2)^2}{u}$, obtained by the quantum Hamiltonian reduction from the boundary admissible $\hat{sl}_2$-modules from Example \ref{ex1}. 
For $j=u-1$ one gets 0, for $u=3$ and $j=0,1$ one gets the trivial representation, but for all other $j$ and $u\geq 5$ the characters are the product sides of the Gordon generalizations of the Rogers-Ramanujan idenities 
(the latter correspond to $u=5$). Another example is 
the minimal series representations of the $N=2$ superconformal algebras,
see \cite{KRW03}, Section 7.
\end{remark}

\end{document}